\documentstyle[12pt, leqno]{article}
\title{Quantization formula for singular reductions}

\author{Youliang Tian\thanks{Partially supported by a NSF postdoctoral
fellowship and a NYU research challenge fund grant.}
$\ $ and Weiping Zhang\thanks{Partially supported by  the NNSF,
SEC of China and the Qiu Shi Foundation.}}

\date{June, 1997}
\begin{document}
\maketitle

\begin{abstract}
In this note, we prove a quantization formula for singular
reductions. The main result is obtained as a simple application of an 
extended quantization formula proved in [TZ2].
\end{abstract}

{\bf \S 0. Introduction}

$\ $

Let $(M,\omega)$ be a closed symplectic manifold. We make the assumption 
that there is a Hermitian line bundle 
$L$ over $M$ admitting a Hermitian connection $\nabla ^L$ with the property that
${\sqrt{-1}\over 2\pi }(\nabla ^L)^2 =\omega$. When such a line bundle exists, it
is usually called a prequantum line bundle over $M$.
Let $J$ be an almost complex
structure on $TM$ so that $g^{TM}(u,v)
=\omega (u,Jv)$ defines a Riemannian metric on $TM$.
Then one can construct canonically a Spin$^c$-Dirac operator
$$D^L: \Omega ^{0,*}(M,L)\rightarrow \Omega ^{0,*}(M,L),
\eqno(0.1)$$
which gives the finite dimensional virtual vector space
$$Q(M,L)=(\ker D^L)\cap \Omega ^{0,{\rm even}}(M,L)-(\ker D^L)
\cap \Omega ^{0,{\rm odd}}(M,L).\eqno (0.2)$$

Next, suppose that $(M,\omega)$ admits a Hamiltonian action of
a compact connected Lie group $G$ with Lie algebra ${\bf g}$.
Let $\mu: M\rightarrow {\bf g}^*$ be the corresponding moment map.
Then a formula due to Kostant [Ko] (cf. [TZ2, (1.13)]) induces a natural
${\bf g}$ action on $L$. We make the assumption that this ${\bf g}$
action can be lifted to a $G$ action on $L$. Then this $G$ action
preserves $\nabla^L$. One can also assume, 
after an integration over $G$ if necessary, that $G$ preserves the Hermitian 
metric on $L$, the almost complex structure $J$ and thus also the 
Riemannian metric $g^{TM}$. $Q(M,L)$ then
becomes a virtual representation of $G$. We denote by
$Q(M,L)^G$ its $G$-invariant part.

Now let $a\in {\bf g}^*$ be a regular value of $\mu$. Let
${\cal O}_a\subset {\bf g}^*$ be the coadjoint orbit of $a$.
For simplicity,
we assume that $G$ acts on $\mu^{-1}({\cal O}_a)$ freely. Then the quotient 
space $M_{G,a}=\mu^{-1}({\cal O}_a)/G$ is smooth. 
Furthermore, $\omega$ descends canonically
to a symplectic form $\omega_a$ on $M_{G,a}$ so that one gets the
Marsden-Weinstein reduction $(M_{G,a},\omega_a)$. The pair $(L,\nabla^L)$
also descends canonically to a Hermitian line bundle $L_{G,a}$ with a
Hermitian  connection  denoted by $\nabla^{L_{G,a}}$. 
\footnote{However, if $a\neq 0$, then 
$L_{G,a}$ is in general no longer a prequantum line bundle over
$(M_{G,a}, \omega_{G,a})$.}  The almost complex structure $J$ also
descends canonically to an almost complex structure on $TM_{G,a}$.
Thus once again one can construct canonically
a Spin$^c$-Dirac operator $D^{L_{G,a}}$ as well as the corresponding
virtual vector space $Q(M_{G,a}, L_{G,a})$. 

The purpose of this short note is to give a direct  proof of
the following result as an immediate consequence of a result proved in
[TZ2].

$\ $

{\bf Theorem 0.1}.  {\it If $\mu^{-1}(0)$ is nonempty, then 
there exists an open neighborhood ${\bf O}$ of $0\in {\bf g}^*$ such that
for any regular value  $a\in {\bf O}$ of $\mu$ with $\mu^{-1}(a)$
nonempty,  the following identity holds,
$$\dim Q(M,L)^G=\dim Q(M_{G,a},L_{G,a}).\eqno(0.3)$$}

If $0\in {\bf g}^*$ is a regular value of $\mu$, then one can take
$a=0$ in (0.3). This is the Guillemin-Sternberg geometric
quantization conjecture [GS], which was proved in various generalities in
[DGMW, GS, G, JK, M1, M2, V1, V2, TZ1, TZ2]. Thus the main feature of
Theorem 0.1 is that it equally applies when $0\in {\bf g}^*$ is a
singular value of $\mu$.

A complete treatment of Theorem 0.1
when $G$ is the circle group was given in [TZ3] using
the spectral flow technique. A different treatment of (0.3)
which works for 
nonabelian actions has also been
given  by Meinrenken and Sjamaar [MS].

Theorem 0.1 still holds,  for orbifolds, when $G$ does not act freely
on $\mu^{-1}({\cal O}_a)$. Furthermore, if $(M,\omega)$ is K\"ahler and
$G$ acts on $M$ holomorphically, then we can refine (0.3) to a system
of holomorphic Morse type inequalities as in [TZ1, 2].

This note is organized as follows. In Section 1, we recall
the result in [TZ2] to be used in this paper. In Section 2, we prove
a key estimate which allows one to obtain Theorem 0.1 immediately.

$$\ $$

{\bf \S 1. An extended quantization formula}

$\ $

We make the same assumption and use the same notation as in Introduction.
We further assume in this section that $0\in {\bf g}^*$ is a regular
value of $\mu$ and, for simplicity, that $G$ acts on $\mu^{-1}(0)$
freely. Then one can construct the canonical Marsden-Weinstein reduction
$(M_G=\mu^{-1}(0)/G,\omega_G)$.

We equip ${\bf g}$ (and thus ${\bf g}^*$ also) 
with an Ad$\, G$-invariant
metric.
Let $h_i$, $1\leq i\leq \dim G$, be an orthonormal base of ${\bf g}^*$. 
Let $V_i$, $1\leq i\leq \dim G$, be the dual base of $h_i$, 
$1\leq i\leq \dim G$.
Then the moment map $\mu$ can be written as 
$$\mu=\sum_{i=1}^{\dim G}\mu_ih_i, \eqno (1.1)$$
with each $\mu_i$ a real function on $M$.

Let 
$${\cal H}=|\mu|^2\eqno (1.2)$$
be the norm square of the moment map. It is clearly $G$-invariant.

For any $V\in {\bf g}$, we use the same notation to denote the vector
field it generates on $M$. 

Let now $E$ be a Hermitian vector bundle over $M$
admitting a Hermitian connection $\nabla ^E$. We make the assumption that
the $G$ action on $M$ lifts to a $G$ action on $E$ which preserves
the Hermitian metric and the Hermitian connection $\nabla ^E$ on $E$.
Then $E$ descends canonically to a Hermitian vector bundle $E_G$ over
$M_G$ with Hermitian connection $\nabla ^{E_G}$.

With these data one can construct canonically a Spin$^c$-Dirac operator
$D^E$ acting on $\Omega^{0,*}(M,E)=\Gamma (\wedge^{0,*}(T^*M)\otimes E)$ 
as well as the corresponding 
virtual vector space $Q(M,E)$ as in (0.2), with $L$ there
replaced now by $E$ (cf. [TZ2, Sect. 1]). Similarly, one can construct
the Spin$^c$-Dirac operator $D^{E_G}$ on $M_G$ as well as the corresponding
virtual vector space $Q(M_G,E_G)$.
Furthermore, one verifies directly that the $G$ action commutes with
$D^E$ so that $Q(M,E)$ is a virtual $G$-representation. Denote its
$G$-invariant part by $Q(M,E)^G$.

Let $V\in {\bf g}$. Set
$$r_V^E=L_V^E-\nabla_V^E,\eqno (1.3)$$
where $L_V^E$ denotes the infinitesimal action of $V$ on $E$.

The following result was proved in Tian-Zhang [TZ2], where
a direct analytic approach to the Guillemin-Sternberg geometric
quantization conjecture [GS] was developed.

$\ $

{\bf Theorem 1.1.} ([TZ2, Theorem 4.2]).  
{\it If $\mu^{-1}(0)$ is nonempty and if at each
critical point $x\in M\backslash \mu^{-1}(0)$ of ${\cal H}$,
$$\sqrt{-1}\sum_{i=1}^{\dim G}\mu_i(x)r_{V_i}^E(x)\geq 0,\eqno (1.4)$$
then the following identity holds,
$$\dim Q(M,E)^G=\dim Q(M_G,E_G).\eqno (1.5)$$}

$\ $

{\bf \S 2. A proof of Theorem 0.1}

$\ $

In this section, we prove Theorem 0.1 by using Theorem 1.1. In order to
explain the basic idea clearer, we will first treat the case where $G$ is
abelian in detail. Then as will be seen, the basic estimate in the 
$G$ nonabelian case, which is necessary for applying Theorem 1.1, can be
reduced to the abelian case.

This section is organized as follows. In a), we deal with the $G$ 
abelian case. In b), we prove Theorem 0.1 for the $G$ nonabelian case.
The final subsection c) contains some further remarks related to
Theorem 0.1.

$\ $

{\bf a) The abelian case}

$\ $

In this subsection, we assume that $G$ is a torus.

For any $a\in {\bf g}^*$, set 
$$\mu_a=\mu-a.\eqno (2.1)$$
Then $\mu_a$ is again a moment map for the $G$-action on $(M,\omega)$.
Let 
$${\cal H}_a=|\mu_a|^2 \eqno (2.2)$$
be the norm square of $\mu_a$.

Recall that $h_i$, $1\leq i\leq \dim G$, is an orthonormal base of 
${\bf g}^*$. Thus $a$ has the expression
$$a=\sum_{i=1}^{\dim G}a_ih_i.\eqno (2.3)$$

The main result of this subsection can be stated as follows.

$\ $

{\bf Proposition 2.1}. {\it There is an open neighborhood ${\bf O}$ of 
$0\in {\bf g}^*$
such that for any $a\in {\bf O}$, 
the following inequality holds at each critical point 
of ${\cal H}_a$,
$$\sum_{i=1}^{\dim G}(\mu_i-a_i)\mu_i\geq 0.\eqno (2.4)$$}

{\it Proof}. We will use Kirwan's geometric
characterization [K1] of the critical
point set of the norm square of moment maps to prove (2.4).

As in [K1, 3.4], denote by ${\bf A}\subset {\bf g}^*$ 
the finite set of weights associated to $\mu$, which consists of
the values of $\mu$ taking on the fixed point set of $G$ on $M$.
We also denote by ${\bf Y}$ the set of convex hulls in ${\bf g}^*$
generated by  nonempty subsets of ${\bf A}$. Then ${\bf Y}$ consists 
naturally of two parts: the part ${\bf Y}_I$ consists of all
those convex hulls not containing $0$, and the rest part 
denoted by ${\bf Y}_{II}$.

Now let $U_\delta$ be an open ball in ${\bf g}^*$ with center $0$ 
and radius $\delta >0$ such that 
the closure $\bar{U}_\delta$ does not intersect with any convex hull in
${\bf Y}_I$. The existence of $U_\delta$ is clear. 
Set ${\bf O}=U_{\delta /2}$

Let $a\in {\bf O}$. Let ${\bf A}_a=\{ A-a:A\in {\bf A}\}$ be the finite set
of weights associated to $\mu_a$, and ${\bf Y}_a=\{ Y-a:Y\in {\bf Y}\}$ the
associated set of convex hulls. Then ${\bf Y}_a$ consists of two parts
accordingly: ${\bf Y}_{I,a}=\{ Y-a:Y\in {\bf Y}_I\}$ and
${\bf Y}_{II,a}=\{ Y-a:Y\in {\bf Y}_{II}\}$. One verifies easily 
that the closure $\bar{\bf O}$ does not intersect with any convex hull
in ${\bf Y}_{I,a}$.

Let $B_a$ be the open ball 
$B_a=\{ y\in {\bf g}^*:|y+{a\over 2}|<|{a\over 2}|\}$. Clearly,
$B_a\subset {\bf O}$. Thus $B_a$ does not intersect with any convex hull
in ${\bf Y}_{I,a}$.

Now take $Y\in {\bf Y}_{II,a}$. Let $y\in Y$ be the (unique) point on $Y$
which is closest to $0$. We claim that 
$y\in {\bf g}^*\backslash B_a$.

To prove this claim, we suppose on the contrary that $y\in B_a$. Let 
$x\in \partial B_a$, with $x\neq -a$,
lie in the straight line generated by $y$ and
$-a\in \partial B_a\cap Y$, then it is easy to see that $x\in Y$. 
For if $x\notin Y$, then $y$ should lie in a face of $Y$ which does not
contain $-a$. This would imply
that $y$ lies in a convex hull in ${\bf Y}_{I,a}$,
a contradiction. But with such an $x\in \partial B_a\cap Y$, one encounters
another contradiction  $d(0,x)<d(0,y)$.
Thus we should have $y\notin B_a$.

By this and by the result of Kirwan [K1, 3.12], one finds that if
$y$ is a critical point of ${\cal H}_a$, then
$$0\leq |\mu_a(y)+{a\over 2}|^2-|{a\over 2}|^2=\sum_{i=1}^{\dim G}
(\mu_i(y)-a_i)\mu_i(y),\eqno (2.5)$$
which is exactly (2.4). $\Box$

$\ $

We can now prove Theorem 0.1 in this abelian case. In fact, 
the Kostant formula [Ko] (cf. [TZ2, (1.13)]) implies, in using the notation
of (1.3), that for each $1\leq i\leq \dim G$,
$$r^L_{V_i}=-2\pi \sqrt{-1}\mu_i.\eqno (2.6)$$
Thus one verifies via Proposition 2.1 that at any critical point
of ${\cal H}_a$, one has
$$\sqrt{-1}\sum_{i=1}^{\dim G}(\mu_i-a_i)r^L_{V_i}
=2\pi \sum_{i=1}^{\dim G}(\mu_i-a_i)\mu_i\geq 0.\eqno (2.7)$$
Theorem 0.1 then follows by applying Theorem 1.1 to 
$(M,\omega,\mu_a,L)$. $\Box$

$$\ $$

{\bf b) The nonabelian case}

$\ $

In this subsection, we no longer assume that $G$ is abelian.

Let  $T$ be a maximal torus of $G$, with Lie algebra ${\bf t}$. Then
$\mu_T=P_T\mu$, where $P_T$ is the orthogonal projection from ${\bf g}^*$
to ${\bf t}^*$, is the moment map of the induced $T$ action on $(M,\omega)$
(cf. [K1, 3.3]). Let ${\bf O}_T\subset {\bf t}^*$ be the open set
defined in a) for  $\mu_T$.  Then 
${\bf O}={\rm Ad}\, G ({\bf O}_T)\subset {\bf g}^*$ is an open 
neighborhood of $0\in {\bf g}^*$.

Let $a\in {\bf O}$ be a regular value of $\mu$ and, for simplicity,
assume that $G$ acts on $\mu^{-1}({\cal O}_a)$ freely.
Recall that the coadjoint orbit ${\cal O}_a$ admits a 
canonical symplectic (actually K\"ahler)
form $\omega_a$ and the Ad$\, G$ action on 
${\cal O}_a$ is Hamiltonian with the moment map given by the 
canonical embedding $i_a:{\cal O}_a\hookrightarrow {\bf g}^*$ (cf. 
[McS, Chap. 5]).

We  now form the symplectic product 
$(M\times {\cal O}_a,\omega\times (-\omega_a))$. Then the induced action
of $G$ on $M\times {\cal O}_a$ is Hamiltonian with the moment map
$\hat{\mu}:M\times {\cal O}_a\rightarrow {\bf g}^*$ given by
$$\hat{\mu}(x,b)=\mu(x)-b.\eqno (2.8)$$
Then $0\in {\bf g}^*$ is a regular value of $\hat{\mu}$ and 
$G$ acts on $\hat{\mu}^{-1}(0)$ freely. Furthermore,
one has the standard identification of the symplectic 
quotients 
$$\hat{\mu}^{-1}(0)/G \equiv \mu^{-1}({\cal O}_a)/G
=M_{G,a}\eqno (2.9)$$
(cf. [McS, Chap. 5]).

Let 
$$\hat{\cal H}=|\hat{\mu}|^2\eqno (2.10)$$
be the norm square of the moment map $\hat{\mu}$.

We now state the nonabelian extension of Proposition 2.1 as follows.

$\ $

{\bf Proposition 2.2}. {\it If $(x,b)\in M\times {\cal O}_a$
is a critical point of $\hat{\cal H}$, then one has the following 
inequality for the inner product on ${\bf g}^*$,
$$\langle \mu(x)-b,\mu(x)\rangle \geq 0.\eqno (2.11)$$}

{\it Proof}. Without loss of generality, we can assume that 
$a\in {\bf O}_T$. 
By a result of Kirwan [K2, pp. 551], we know that for any
critical point
$(x,b)\in M\times {\cal O}_a$
of $\hat{\cal H}$ , one
can find $y\in {\bf t}^*$ such that $(y,a)$ is a critical
point of     $\hat{\cal H}$ in the $G$-orbit of $(x,b)$.

Now since $a\in {\bf O}_T$,
one finds from Proposition 2.1 and from Kirwan [K1, 3.3]
that
$$\langle \mu (y)-a, \mu (y)\rangle \geq 0.\eqno (2.12)$$
(2.11) then follows from the Ad$\, G$-invariance of the inner product
on ${\bf g}^*$. $\Box$

$\ $

{\it Proof of Theorem 0.1}. Denote by $\pi$ the projection from
$M\times {\cal O}_a$ to its first factor $M$. Let
${\cal L}=\pi^*L$ be the pull-back Hermitian
line bundle over  $M\times {\cal O}_a$
with the pull-back Hermitian connection $\nabla^L$
on ${\cal L}$ verifying  that
$${\sqrt{-1}\over 2\pi}(\nabla^{\cal L})^2=\pi^*\omega. \eqno (2.13)$$
Furthermore, the $G$ action on $L$ lifts canonically
to an action on ${\cal L}$. In particular, for any $V\in {\bf g}$,
its induced infinitesimal action on ${\cal L}$ is, via
the Kostant formula [Ko] (cf. [TZ2, (1.13)]) for the ${\bf g}$-action
on $L$,    given by
$$L^{\cal L}_V=\nabla^{\cal L}_V-2\pi \sqrt{-1}\langle \mu\pi,
V\rangle ,\eqno (2.14)$$
from which one has, in using the notation  in (1.3), that
$$r^{\cal L}_{V_i}=-2\pi \sqrt{-1} \mu_i(x)\eqno (2.15)$$
at any point $(x,b)\in M\times {\cal O}_a$.
Thus by (2.15) and Proposition 2.2, one verifies that at any critical point
$(x,b)\in M\times {\cal O}_a$
of $\hat{\cal H}$, 
$$\sum_{i=1}^{\dim G} \sqrt{-1}\hat{\mu}_i(x,b)r^{\cal L}_{V_i}
(x,b)= 2\pi \langle
\mu(x)-b,\mu(x)\rangle \geq 0.\eqno (2.16)$$

One the other hand, one verifies directly that the induced
line bundle ${\cal L}_G$ over $\hat{\mu}^{-1}(0)/G$ is exactly the line
bundle $L_{G,a}$   over $M_{G,a}=\mu^{-1}({\cal O}_a)/G$.

One can then apply  Theorem 1.1 to $M\times {\cal O}_a$,
$\hat{\mu}$ and ${\cal L}$ to get
$$\dim Q(M\times {\cal O}_a, {\cal L})^G
=\dim Q(M_{G,a},L_{G,a}).\eqno (2.17)$$
Furthermore, by the definition of ${\cal L}$, one verifies directly
that
$$\dim Q(M\times {\cal O}_a,{\cal L})^G=\dim Q(M,L)^G \cdot
\dim Q({\cal O}_a,{\bf C})^G, \eqno (2.18)$$
with
$$\dim Q({\cal O}_a, {\bf C})^G=1,\eqno (2.19)$$
which can also be verified directly.
(0.3)  follows from (2.17)-(2.19).  $\Box$

$\ $

{\bf c) Further extensions and remarks}

$\ $

First of all, as have already been remarked in Introduction,
Theorem 0.1 still holds for orbifolds when $G$ does not
act on $\mu^{-1}({\cal O}_a)$ freely, and that in the holomorphic
category one can refine (0.3) to a system of Morse type
inequalities. We will not fill the easy details here.

The second point we want to remark is that by now it is clear
that Theorem 1.1 has a very wide range of applicability, and
Theorem 0.1 is just one manifestation of this. One can easily
state and prove an extended version of Theorem 0.1 which works for
general coefficients. In particular, one can take the trivial
line bundle  as coefficient to get
as in [TZ2] that the Todd genus of $M$ equals to the Todd genus
of $M_{G,a}$ for any regular value $a$ of $\mu$. This
last result is also contained in [MS].

Last but not least, recall that we have extended the
Guillemin-Sternberg geometric quantization conjecture [GS]
to the case of manifolds with boundary in [TZ3].
Combining the methods there
with our proof of Theorem 0.1, one can easily extend
Theorem 0.1 to the case of manifolds with boundary as well. 
We again leave the details
to the interested reader.

$$\ $$

{\bf References}

$\ $

[DGMW] H. Duistermaart, V. Guillemin, E. Meinrenken and S. Wu,
Symplectic reduction and Riemann-Roch for circle actions.
{\it Math. Res. Lett}.
2 (1995), 259-266.

[G] V. Guillemin, Reduced phase space and Riemann-Roch. in
{\it Lie Groups and
Geometry}, R. Brylinski et al. ed., Birkhaeuser
(Progress in Math., 123),
1995, pp. 305-334.

[GS] V. Guillemin and S. Sternberg, Geometric quantization and
multiplicities of group representations.
{\it Invent. Math.} 67 (1982), 515-538.

[JK] L. C. Jeffrey and F. C. Kirwan, Localization and quantization
conjecture. {\it Topology} 36 (1997), 647-693.

[K1] F. C. Kirwan, {\it Cohomology of Quotients in Symplectic and Algebraic
Geometry.}  Princeton Univ. Press, Princeton, 1984.

[K2] F. C. Kirwan, Convexity properties of the moment mapping, III.
{\it Invent. Math.} 77 (1984), 547-552.

[Ko] B. Kostant, Quantization and unitary representations. {in} {\it Modern
Analysis and Applications.}  Lecture Notes in Math. vol. 170,
Springer-Verlag, (1970), pp. 87-207.

[M1] E. Meinrenken, On Riemann-Roch formulas for multiplicities.
{\it J. Amer. Math. Soc.}
9 (1996), 373-389.

[M2] E. Meinrenken, Symplectic surgery and the Spin$^c$-Dirac operator.
{\it  Adv. in Math.} To appear.

[MS] E. Meinrenken and R. Sjamaar, Singular reduction and quantization.
{\it Preprint,} 1996.

[McS] D. Mcduff and D. Salamon, {\it Introduction to Symplectic
Topology}. Clarendon Press, Oxford, 1995.

[TZ1] Y. Tian and W. Zhang, Symplectic reduction and quantization.
{\it C. R. Acad. Sci. Paris,  S\'erie I}, 324 (1997), 433-438.

[TZ2] Y. Tian and W. Zhang, An analytic proof of the geometric
quantization conjecture of Guillemin-Sternberg.
{\it Preprint}, 1997. (This is a revised version of
``Symplectic reduction and analytic localization.
{\it Preprint}, 1996".)

[TZ3] Y. Tian and W. Zhang, Quantization formula for symplectic manifolds
with boundary. {\it Preprint}, 1997.

[V1] M. Vergne, Multiplicity formula for geometric quantization, Part I.
{\it Duke Math. J}. 82 (1996), 143-179.

[V2] M. Vergne, Multiplicity formula for geometric quantization, Part II.
{\it Duke Math. J}. 82 (1996), 181-194.

$\ $

Y. T.: Courant Institute of Mathematical Sciences,
New York, NY 10012,
U. S. A.

{\it e-mail address:}
ytian@cims.nyu.edu

$\ $

W. Z.:
Nankai Institute of Mathematics,
Tianjin, 300071,
P. R. China

{\it e-mail address:}
weiping@sun.nankai.edu.cn
\end{document}